\documentclass[reqno]{amsart}
\usepackage[margin=1in]{geometry}

\usepackage[english]{babel}
\usepackage{a4wide,xcolor,graphicx,amsmath,amssymb,verbatim,mathrsfs,latexsym}
\usepackage[colorlinks=true, pdfstartview=FitV,linkcolor=blue,citecolor=blue, urlcolor=black]{hyperref}
\usepackage{subcaption} 
\usepackage{float}      
\usepackage{cleveref}
\usepackage{dsfont}
\usepackage{booktabs,multirow,array,makecell}
\usepackage{enumitem}
\crefformat{equation}{(#2#1#3)}
\usepackage{placeins}

\crefname{section}{Section}{Sections}
\Crefname{section}{Section}{Sections}

\crefname{theorem}{Theorem}{Theorems}
\Crefname{theorem}{Theorem}{Theorems}

\crefname{lemma}{Lemma}{Lemmas}
\Crefname{lemma}{Lemma}{Lemmas}
  
\crefname{proposition}{Proposition}{Propositions}
\Crefname{proposition}{Proposition}{Propositions}

\crefname{corollary}{Corollary}{Corollaries}
\Crefname{corollary}{Corollary}{Corollaries}

\crefname{example}{Example}{Examples}
\Crefname{example}{Example}{Examples}

\crefname{definition}{Definition}{Definitions}
\Crefname{definition}{Definition}{Definitions}

\crefname{remark}{Remark}{Remarks}
\Crefname{remark}{Remark}{Remarks}

\crefname{proposition*}{Proposition}{Propositions}
\Crefname{proposition*}{Proposition}{Propositions}
\crefname{remark*}{Remark}{Remarks}
\Crefname{remark*}{Remark}{Remarks}
\crefname{claim*}{Claim}{Claims}
\Crefname{claim*}{Claim}{Claims}

\allowdisplaybreaks
\numberwithin{equation}{section}
\newtheorem{theorem}{Theorem}[section]
\newtheorem{proposition}{Proposition}[section]
\newtheorem{lemma}{Lemma}[section]

\newtheorem{corollary}{Corollary}[section]
\newtheorem{remark}{Remark}[section]

\newcommand{\init}{\textnormal{in}}

\title[Homoenergetic Landau equation]{The Landau equation: Well-posedness and singularities of homoenergetic solutions}

\author{Francisco-Unai Caja-Lopez}
\address{The University of Texas at Austin, USA}
\email{funai.caja@utexas.edu}

\author{Maria Pia Gualdani}
\address{The University of Texas at Austin, USA}
\email{gualdani@math.utexas.edu}

\author{Alessia Nota}
\address{Gran Sasso Science Institute, L’Aquila, Italy}
\email{alessia.nota@gssi.it}

\thanks{MPG is partially supported by NSF grant DMS-2511625 and by the CNS SPARK grant.}
\date{\today}

\begin{document}


\begin{abstract}

We consider the inhomogeneous Landau equation and focus on solutions whose density and temperature evolve in time but remain spatially homogeneous. Such solutions are called homoenergetic. We identify one class of initial data for which homoenergetic solutions exist globally in time, and another class for which they develop finite-time hydrodynamic singularities. The initial data are spatially inhomogeneous, neither close to vacuum nor close to equilibrium; both the well-posedness and the finite-time blow-up results hold for all soft and hard potentials.


\end{abstract}

\maketitle
 \setcounter{tocdepth}{2}
\tableofcontents

\section{Introduction}

The Landau equation is a fundamental model in kinetic theory which can be used to describe the dynamics of collisional plasmas. It characterizes the evolution of a large system of interacting particles through the evolution of the phase-space density $f(t,x,v)$, where $t \geq 0$ denotes time and $x,v \in \mathbb{R}^3$ represent the position and velocity variables, respectively.
The Landau equation is given by
\begin{equation}\label{eq:L}
\partial_t f+v\cdot\nabla_x f=Q(f),
\end{equation}
where the collision operator acts only on the velocity variable and is defined by
\begin{equation}\label{eq:collopL}
Q(f)(v)=\operatorname{div}_v\int_{\mathbb R^3}\Phi(v-v_*)
\Big(f(t,x,v_*)\nabla_{v} f(t,x,v)-f(t,x,v)\nabla_{v_*}f(t,x,v_*)\Big)\,dv_*,
\end{equation}
with
\begin{equation}
    \label{eq:opPhi}
\Phi(z)=|z|^{2+\gamma}\mathbb P(z),\qquad
\mathbb P(z)=\operatorname{Id}-\frac{z\otimes z}{|z|^2},
\end{equation}
where $\mathbb P(z)$ denotes the orthogonal projection onto the plane orthogonal to $z$. The parameter $\gamma$  is usually considered in the range $\gamma\ge -3$, with the most physically relevant case corresponding to $\gamma=-3$. This case is usually referred to as Coulomb interaction. 

The question of whether solutions to \eqref{eq:L} develop finite-time singularities remains a challenging open problem. Global-in-time smooth solutions are known to exist only for specific  classes of initial data. In the spatially homogeneous setting, this was established by Guillen and Silvestre \cite{GS24}, with subsequent developments for non smooth initial data by Golding, Gualdani, and Loher \cite{GualdGoldLoher2025LandauL32}, Ji \cite{Ji2024DisFisher}, and Desvillettes et al. \cite{GualdaGoldLohDev2024LandauFisherL1}. In the perturbative regime near equilibrium, global well-posedness was proved by Guo \cite{Guo2002LandauPeriodicBox} and later refined by Guo and Strain \cite{Strain2008exponential}. For initial data close to vacuum and interaction potentials satisfying $\gamma>-2$, global existence was established by Luk \cite{luk2019vacuumsoft} for soft potentials and by Chaturvedi \cite{chaturvedi2022vacuumhard} for hard potentials. On the other hand, the second author et al. \cite{bedrossian2022non} ruled out a class of Type I self-similar blow-up profiles. More recently, the second author et al. \cite{Gualdani2026Blowup} proved the existence of finite-time singularities for interaction potentials $\gamma \in (\sqrt{3},2]$. Their proof exploits a connection between the Landau equation and the compressible isentropic Euler equations. 

The contribution of this manuscript is twofold: we show the existence of a class of smooth initial data for which the equation (\ref{eq:L}) admits a global-in-time smooth solution, and another class for which smooth solutions develop hydrodynamic singularities in finite time. Our results are collected in the following two theorems:

\begin{theorem}[Global well-posedness]\label{thm:main}
Let $T>0,\; \gamma \geq -3$. Let $L_0\in \mathbb{R}^{3\times 3}$ satisfy 
\begin{align}\label{cond:A}
\det(I+tL_0)> 0 
\end{align}
for all $0\leq t\leq T$. Let $g_\init(\cdot) \ge 0$ be a smooth and rapidly decaying function in $\mathbb{R}^3$.  

The inhomogeneous Landau equation \cref{eq:L} with initial data  
$$f_\init(x,v)=g_\init(v-L_0x),$$
    has a classical solution ${f\in \mathcal C^\infty\left([0,T]\times\mathbb R^6\right)} \cap L_{t,x}^\infty(L^p_v)$ for all $1\le p\le +\infty$.

\end{theorem}

Our second result concerns formation of hydrodynamic singularities in finite time in the case $L_0=-I$.

\begin{theorem}[Hydrodynamic implosion]\label{thm:blowup_general}
Let $\gamma \geq -3$, and let $g_\init\ge 0$ be a smooth, rapidly decaying function on $\mathbb R^3$ with $\int_{\mathbb R^3} w\,g_\init(w)\,dw=0$. Set
    $$
    f_\init(x,v):=g_\init(x+v).
    $$
    Then the inhomogeneous Landau equation \cref{eq:L} with initial data $f_\init$
    admits a solution $f$ with the following properties:
    \begin{enumerate}[label=\roman*)]
        \item $f$ is smooth and bounded up to time $t=1$: $f\in \mathcal C^\infty\left([0,1)\times \mathbb R^6\right)\cap L^\infty([0,1)\times \mathbb{R}^6)$. 
        \item The hydrodynamic quantities of $f$ blow up at $t=1$. More precisely, the density $\rho$, bulk velocity $u$, and temperature $\theta$ of $f$ are given by
            \begin{align}\label{implos_rho}
                \rho(t,x) := & \int _{\mathbb R^3} f(t,x,v)\,dv = (1-t)^{-3}\int_{\mathbb R^3} g_\init(w)\,dw,\\
                \label{implos_u}
                u(t,x) := & \frac{1}{\rho(t,x)} \int_{\mathbb R^3} v f(t,x,v)\,dv = -(1-t)^{-1}x,\\
                \theta(t,x) := & \frac{1}{3}\frac{\int_{\mathbb R^3}|v-u(t,x)|^2 f(t,x,v)\,dv}{\rho(t,x)} 
                =\frac{1}{3} (1-t)^{-2}\frac{\int_{\mathbb R^3}|w|^2g_\init(w)\,dw}{\int_{\mathbb R^3}g_\init(w)\,dw}.\label{implos_theta}
            \end{align}
The hydrodynamic fields (\ref{implos_rho})-(\ref{implos_u}) are a solution to the monoatomic compressible Euler equations (\ref{comp_mon_Euler}). 

        \item If $\gamma\geq -2$, as $t\to1^{-}$ the solution approaches the time-dependent Maxwellian
        $$
        f(t,x,v) \sim 
        \frac{\rho(t)}{(2\pi\theta(t))^{\frac{3}{2}}}\exp\left\{-\frac{|v-u(t,x)|^2}{2\theta(t)}\right\}.
        $$
        If $\gamma<-2$, as $t\to1^{-}$  the profile approaches an intermediate state between $g_{\init}$ and a Maxwellian: 
        $$
        f(t,x,v) \sim F(s_0, x+(1-t)v), \qquad s_0=-\frac{1}{\gamma+2},
        $$
        where $F(s_0,w)$ is the solution to the homogeneous Landau equation with initial data $g_{\init}(w)$ computed at $s=s_0$. 
    \end{enumerate}
\end{theorem}
\smallskip


All statements made in Theorem \ref{thm:main} and Theorem \ref{thm:blowup_general}  have a clear analog for  the Boltzmann equation. In fact, already in 1963 Nikolskii observed in \cite{Nikolskii1963} that certain solutions to the Boltzmann equation exhibits a compressive behavior similar to the one of \cref{thm:blowup_general}. We have chosen to focus on the Landau equation for the sake of simplicity. 

Some more comments are in order:

{\bf{Initial data.}}  The initial data $f_{\init}$ in both Theorem \ref{thm:main} and Theorem \ref{thm:blowup_general} are smooth and bounded, but not integrable in $\mathbb{R}^3 \times \mathbb{R}^3$. They belong to $L^\infty_xL^1_v(\mathbb{R}^3)\cap L^\infty(\mathbb{R}^6)$. The initial hydrodynamic quantities $\rho_{\init}(x)$ and $\theta_{\init}(x)$ do not depend on $x$; they will evolve in time but will remain homogeneous with respect to $x$.  
The first and second moments of $f_{\init}$ grow in $x$ linearly and quadratically respectively: 

\begin{align*}
\rho_{\init} u_{\init} & = \int v \;f_{\init}(x,v) \;dv = c_1+ \rho_{\init} L_0x, \\
{M_2}(f_{\init})& =\int |v|^2 \; f_{\init}(x,v) \;dv = c_2 + |L_0x|^2 \rho_{\init} + 2 L_0x \cdot c_1,
\end{align*}
where
$$
c_1 := \int w \; g_{\init}(w)\;dw,\qquad c_2 := \int |w|^2 \;g_{\init}\;dw,
$$
and we remind that $L_0\in \mathbb{R}^{3\times 3}$.
We will choose initial data such that $\int_{\mathbb R^3} wg_\init(w)\,dw=0$ to simplify certain computations, but the condition is not necessary to prove the theorem. We stress that we are not in a perturbative setting, we have no restriction on the size of $g_{\init}$ and $L_0$.

{\bf{The main results.}}
The key observation in the proof of  Theorem \ref{thm:main}  is that  solutions to (\ref{eq:L}) with initial data of the form $g_{\init}(v-L_0x)$ have  the form
    \begin{equation}\label{eq:homoenergetic}
        f(t,x,v)=g(t,v-L(t)x), \qquad L(t)=L_0(I+tL_0)^{-1},
    \end{equation}
    and that $g$ satisfies a modified homogeneous equation (see \cref{eq:drifted_homo_Landau} below). Hence, analyzing (\ref{eq:L}) with initial data of the form $g_{\init}(v-L_0x)$ is equivalent to  analyzing  \cref{eq:drifted_homo_Landau} with general initial data. Our solution $f$ will have bounded mass and entropy but {\em{unbounded}} second and higher moments, see Corollary \ref{cor:moments-f}.  Although \cref{thm:main} is stated for smooth initial data, we expect the same conclusion to hold for $g_\init\in L^1$; this should follow by adapting the argument of \cite{Ji2024DisFisher} or \cite{GualdaGoldLohDev2024LandauFisherL1} for the homogeneous Landau equation. Examples of matrices $L_0$ that satisfy the non-degeneracy condition (\ref{cond:A}) globally in time are $L_0 = I$ (isotropic expansion/dilation) or any $L_0 = \operatorname{diag}(\lambda_1, \lambda_2, \lambda_3)$, $\lambda_i \geq 0$ (expansion, possibly anisotropic) or a shear matrix 
    \begin{equation}\label{eq:shear_mat}
        L_0 = \begin{pmatrix} 0 & K & 0 \\ 0 & 0 & 0 \\ 0 & 0 & 0 \end{pmatrix}, \qquad K\in \mathbb{R}_{+}\,.
    \end{equation}
    More examples are given in Appendix \ref{sec:asymptotics}.\\

\begin{figure}[h]
  \centering
  \begin{subfigure}[b]{0.47\textwidth}
    \centering
    \includegraphics[width=\textwidth]{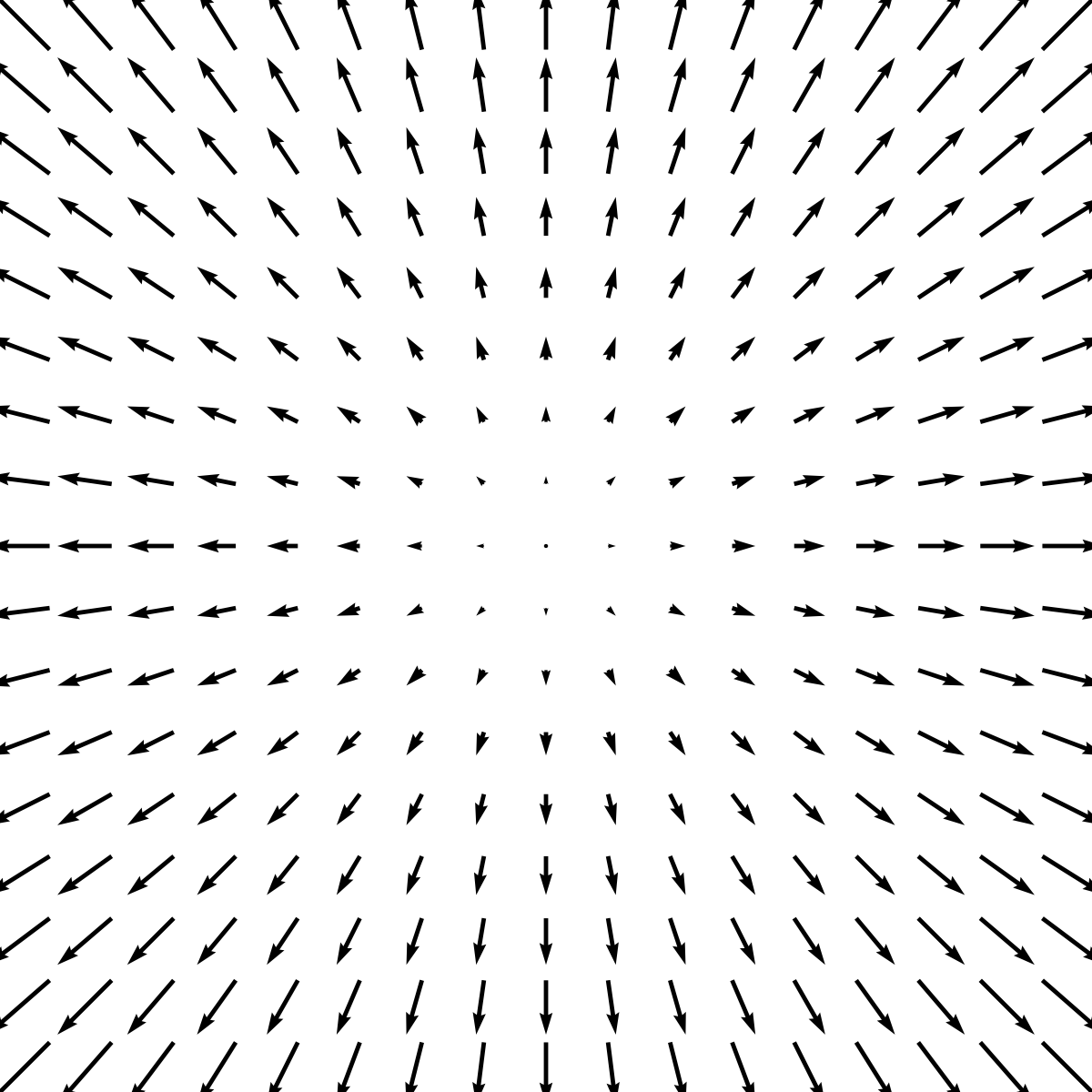}
    \caption{Case with $L_0=I$}
    \label{subfig:minus_v_field}
  \end{subfigure}
  \hfill
  \begin{subfigure}[b]{0.47\textwidth}
    \centering
    \includegraphics[width=\textwidth]{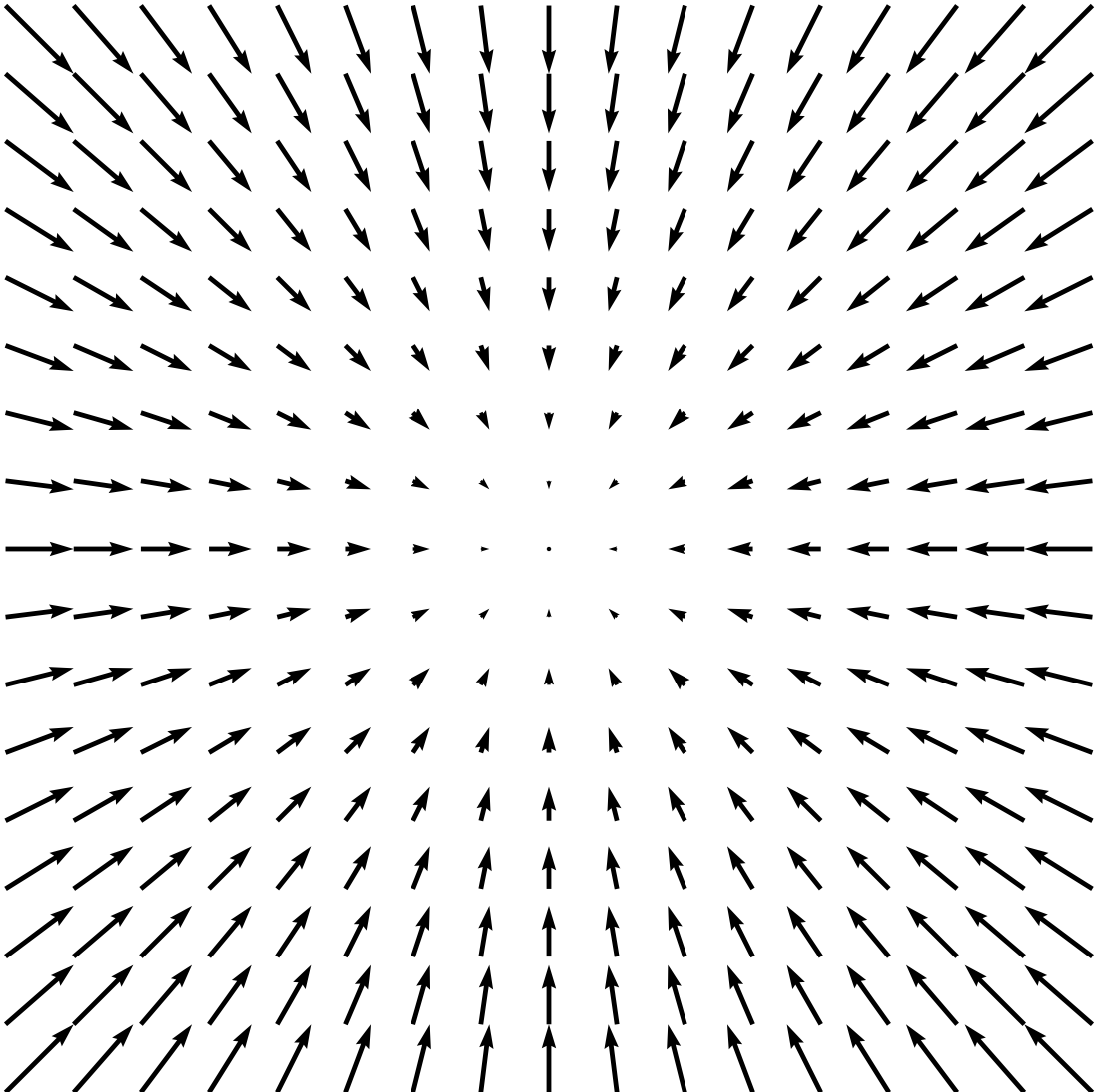}
    \caption{Case with $L_0=-I$}
    \label{subfig:plus_v_field}
  \end{subfigure}
  \hfill
  \caption{Initial bulk velocities $u_0(x) = \pm x$ corresponding to initial data $f_\init(x,v)=g_\init({v}-L_0x)$ with $L_0=\pm I$.}
\end{figure}

The proof of Theorem \ref{thm:blowup_general} has a particularly simple structure. We will show that  $f$ has a semi-explicit expression
    \begin{equation}\label{eq:ansatz_blow_up}
        f(t,x,v)=h\left(s(t),(1-t)v+x\right),\qquad
    s(t):=\int_0^t(1-\tau)^{-3-\gamma}\,d\tau,
    \end{equation}
    where $h$ is the solution to the homogeneous Landau equation with initial data $g_\init$. The structure of the initial data is of the same type as the ones of Theorem \ref{thm:main} if one takes $L_0= -I$ (isotropic compression). For such choice of $L_0$ the condition $\textrm{det}(I+tL_0)>0$ fails at time $t=1$. Notice that  the choice $L_0=I$ falls into the case of Theorem \ref{thm:main}, and solutions to (\ref{eq:L}) with initial data $f_\init(x,v)=g_\init(v-x)$ exhibit drastically different behavior from the ones with $f_\init(x,v)=g_\init(v+x)$. One way to see this qualitative difference is through the bulk velocity 
 $$
 u(t,x):= \frac{\int vf(t,x,v)\,dv}{\int f(t,x,v)\,dv}.
 $$
 For $L_0=I$ at time $t=0$ one has $u(0,x)=x$, whereas for $L_0=-I$ one has $u(0,x)=-x$. Two-dimensional versions of these vector fields are displayed in Figures \ref{subfig:minus_v_field} and \ref{subfig:plus_v_field}, respectively. 

 For $\gamma=-3$, the solution $f$ from Theorem  \ref{thm:blowup_general} actually stays regular in time until  $t=1$ included, $f\in \mathcal C^\infty([0,1]\times\mathbb{R}^6)$. For $\gamma \neq -3$, $f$ loses regularity at the blow up time due to the rescaling of the kernel. See equation \cref{eq:self_similar_final} for more details. The hydrodynamic fields (\ref{implos_rho})-(\ref{implos_u}) are a solution to the monoatomic compressible Euler equations
\begin{equation}\label{comp_mon_Euler}
    \begin{cases}
        \partial_t\rho + \textnormal{div}(\rho u)=0\\
        \partial_t(\rho u) + \textnormal{div}(\rho(u\otimes u)+pI)=0,\\
        \partial_t\left(\rho\left(\frac{1}{2}|u|^2+\frac{3}{2}\theta\right)\right) +
    \textnormal{div}\left(\left(\rho\left(\frac{1}{2}|u|^2+\frac{3}{2}\theta\right)+p\right)u\right) =0,
    \end{cases}\qquad
    p = \rho\theta,
\end{equation}
where $\theta$ denotes the temperature per unit mass. The initial data is $(\rho,u,\theta)|_{t=0}=(\rho_\init, -x, \theta_\init)$, where $\rho_\init,\theta_\init>0$ are fixed constants.


General deformation matrices $L_0$ satisfying $\det(I+t_*L_0)=0$ cannot be related to the homogeneous problem. Nevertheless, we still observe the formation of finite-time singularities in the hydrodynamic quantities, in particular the mass blows up as
\[
\rho(t)=\frac{\rho(0)}{\det(I+tL_0)}\gtrsim (t_*-t)^{-n_{\lambda}},
\]
where $n_\lambda$ is the multiplicity of the eigenvalue $\lambda = -\frac{1}{t_*}$ of $L_0$.

\textbf{Comparison with conditional regularity results.}
Theorem  \ref{thm:blowup_general}, to the best of our knowledge, provides the first example of a solution to the Landau equation that exhibits finite-time singularities in the range ${\gamma\in[-3,\sqrt 3]}$. 
For $\gamma>\sqrt 3$, implosion singularities were recently constructed in \cite{Gualdani2026Blowup} by exploiting the connection of the Landau equation with the compressible Euler system. 

As in \cite{Gualdani2026Blowup}, also in our case the singularity occurs due to tail fattening, meaning the function $f$ remains bounded but particles accumulate at large velocities. This is surprising, since the mechanism of implosion due to tail fattening for soft potentials was recently ruled out in \cite{GHS25}. There is no contradiction, however, between our work and \cite{GHS25}: the a priori estimate in \cite{GHS25} applies only to initial data with velocity decay of the form
\[
f_\init(x,v)\leq C(1+|v|)^{-m},
\]
for some $m>0$, a condition not satisfied by the initial data considered here. Similarly, the conditional results shown by Snelson and Solomon in \cite{SnelsonSolomon2026} do not apply here because the initial data is not localized as required by their theorem.

\textbf{Homoenergetic solutions.}  Solutions arising from Theorems \ref{thm:main} and \ref{thm:blowup_general} are of {\bf{homoenergetic}} type. Homoenergetic solutions to kinetic equations were first  introduced by Galkin and Truesdell \cite{Galkin1958,Truesdell1980, truesdell1956pressures} in the context of the Boltzmann equation. Their defining feature is that even if the function $f$ depends on $v$ and $x$, its density $\rho(t,x)$ and internal energy $\theta(t,x)$ are  homogeneous in space and only depend on the time $t$. Cercignani  in \cite{Cercignani1989} studied homoenergetic affine flows for Boltzmann with hard-sphere and angular-cutoff interactions and showed that, depending on the initial data, solutions exist globally in time or exhibit an implosion with the density becoming unbounded in finite time.

A systematic mathematical theory for the long time regime of the homoenergetic Boltzmann
equation was subsequently developed in \cite{james2019long,james2019self,james2020long, bobylev2020self, nota2026non, NotaVelazque2023Shear}. These works highlight how the presence of the deformation matrix affects the long time asymptotics.

 Several subsequent works have broadened this framework: 
\cite{kepka2023self,kepka2024longtime} addressed non-cutoff Maxwell molecules and hard potentials; 
\cite{duan2021boltzmann,DuanLiu2025} treated uniform shear flow, deriving sharp regularity and stability results; 
\cite{miele2025homoenergetic, miele2026asymptotic} studied the linear Rayleigh--Boltzmann analogue;
and \cite{carrillo2025inelastic} recently extended the theory to dissipative granular collisions under shear heating.
More details about the long time asymptotics together with the matrix classification are provided in Appendix \ref{sec:asymptotics}.

\textbf{Structure of the paper.}  In \cref{sec:derivation}, we present a motivation of the homoenergetic ansatz and derive the equation \cref{eq:drifted_homo_Landau} satisfied by the function $g$ from \cref{eq:homoenergetic}. In \cref{sec:well_posedness} we show well posedness of solutions to  \cref{eq:drifted_homo_Landau}. In \cref{sec:density_blow_up}, we prove \cref{thm:blowup_general}. Finally, in Appendix \ref{sec:asymptotics}, we recall a classification of deformation matrices which remain regular for all times and provide a heuristic argument that predicts the long time behavior of solutions.

\textbf{Notation.} Throughout, $I$ denotes the $3\times3$ identity matrix and, for a matrix $L\in\mathbb R^{3\times 3}$, $|L|$ denotes its operator norm. We write $\langle w\rangle:=(1+|w|^2)^{1/2}$. Boltzmann's entropy, moments of order $s$, and mass are respectively denoted
$$
\mathcal H(g)=\int g\log g\,dw\,\quad
M_s(g)=\int|w|^sg\,dw,\quad 
\rho(t)=M_0(g(t)).
$$
Fisher information is denoted
\begin{equation}\label{def:fisher}
I(g)=\int\frac{|\nabla g|^2}{g}\,dw=4\int|\nabla\sqrt g|^2\,dw=\int|\nabla\log g|^2g\,dw.
\end{equation}
The Gateaux derivative of a functional $\mathcal E(g)$ at $g$ along the direction $h$ is denoted $\langle \mathcal E'(g),h\rangle$.

\textbf{AI statement.} AI tools were used only to improve the exposition and to proofread the final version of the manuscript.

\section{The homoenergetic flow}\label{sec:derivation}

While the early works of Truesdell and Galkin studied homoenergetic flows at the level of the moment equations, the equation governing $g$ was first written explicitly at the level of the Boltzmann distribution function by Cercignani \cite{Cercignani1989} and later by James, the third author, and Velazquez in \cite{james2019self}. In the same spirit, we begin 
with the classical inhomogeneous Landau (or Boltzmann) equation
$$
\partial_{t}f+v\cdot\nabla_{x}f=Q(f),
$$
and plug in the ansatz
$$
f(t,x,v)=g\left(t, v-\xi(t,x)\right),
$$
which yields the following equation for $g$:
$$
\partial_{t}g-\left( D\xi\cdot w+D\xi\cdot\xi+\partial_{t}\xi\right)\cdot\nabla_{w}g =Q(g).
$$
Here 
$$
w := v-\xi(t,x). 
$$
We are looking for deformations $\xi(t,x)$ such that the quantity $\partial_{t}\xi+D_{x}\xi\cdot w+D_{x}\xi\cdot\xi$ depends on $t$, and on $x$ and $v$ solely through $w$:

\begin{equation}\label{eq:ansatz_deriv}
    \partial_{t}\xi+D_{x}\xi\cdot w+D_{x}\xi\cdot\xi=Z(t,w).
\end{equation}
We show that the only $Z(t,w)$ that satisfies condition (\ref{eq:ansatz_deriv}) is a vector affine in $w$: first note that $Z(t,w)-D_{x}\xi\cdot w$ is independent of $w$, which implies that
$$
D_{w}Z(t,w)=D_{x}\xi(t,x).
$$
Therefore 
$$
Z(t,w)=D_{x}\xi(t,x)w+a(t),\qquad D_{x}\xi(t,x)= L(t),
$$
for some functions $a(t)$  and matrix $L(t)\in\mathbb{R}^{3\times3}$. Hence 
\begin{align}\label{ansatz:xi}
\xi(t,x)=L(t)x+b(t),
\end{align}
and 
\begin{align}\label{ansatz:A}
Z(t,w)= L(t)w+a(t). 
\end{align}
Plugging (\ref{ansatz:xi}) and (\ref{ansatz:A}) into  \cref{eq:ansatz_deriv} yields 
$$
\left(L'(t)+L(t)^2 \right)x+b'(t)+L(t)b(t)=a(t),
$$
which implies  
$$
L'(t) = -L(t)^2,\qquad b'(t)=a(t)  {-L(t)b(t)}.
$$
The equation for $L$ has the explicit solution
$$
L(t)=L_0(I+tL_0)^{-1},
$$
which leaves us with the final simplification
$$
f(t,x,v)=g(t,v-L(t)x-b(t)),\qquad\partial_t g -\left(L(t)w+a(t)\right)\cdot\nabla_wg = Q(g).
$$
Finally, the change of variables $\tilde g(t,w)=g(t,w-b(t))$ reduces the equation to the case $a(t)=b(t)=0$:
\begin{equation}\label{eq:drifted_homo_Landau}
    \partial_{t}g=Q(g)+L(t)w\cdot\nabla_{w}g.
\end{equation}

\section{Well posedness for the homoenergetic flow}\label{sec:well_posedness}

\subsection{A priori estimates and Fisher information}\label{sec:apriori}

This section contains a-priori estimates for solutions to \cref{eq:drifted_homo_Landau} with initial data $g_{\init}(w)$
and $L(t) = L_0(I+tL_0)^{-1}$ for some matrix $L_0$.  The lemmas in this section hold in any time interval $[0,T]$  as long as 
$$
\textrm{det}(I+ tL_0)>0\quad \textrm{ for all}\;  t\in [0,T]
$$
and for any $\gamma\ge-3$.

\begin{lemma}[Moments and Entropy] \label{lem:moments_entropy}
    Let $g$ be a solution of \cref{eq:drifted_homo_Landau} with $M_2(g_\init)<\infty$, $\mathcal H(g_\init)<\infty$. The function $g$ satisfies the following properties: 
    \begin{enumerate}
\item The mass density $\rho(t)$ satisfies the identity
    \begin{equation}\label{eq:mass}
        \rho(t)=\rho(0)\exp\left\{ -\int_{0}^{t}\textnormal{tr}\,L(s)\,ds\right\} = \frac{\rho(0)}{\det (I+tL_0)}, \qquad \rho(0) := \int_{\mathbb{R}^3} g_{\init}\;dw.
    \end{equation}
   \item  The first moment $\vec M_1(g):=\int w\,g(w)\,dw$ satisfies the differential equation
    \begin{equation}\label{eq:M1}
        \frac{d}{dt}\vec M_1(g) = -\textnormal tr\,L(t)\vec M_1(g) - L(t)\vec M_1(g),
    \end{equation}
    In particular, if $\vec M_1(g_\init)=0$ then $\vec M_1(g(t))=0$ for all $t$.
    
    \item The second moment satisfies the estimate
    \begin{equation}\label{eq:second_moment_bound}
        M_{2}(g)\leq M_{2}(g_{\init})\exp\left\{ \int_{0}^{t}\Big(2\left|L(s)\right|-\textnormal tr\,L(s)\Big)\,ds\right\}.
    \end{equation}

    \item 
    The entropy satisfies the identity
    \begin{equation}\label{eq:entropy_bound}
    \mathcal{H}\left(g(t)\right)+\rho(t)\int_{0}^{t}\rho(s)^{-1}D_h\left(g(s)\right)\,ds=\frac{\rho(t)}{\rho(0)}\mathcal{H}(g_{\init}),
\end{equation}
    where $D_h(g)$ is the entropy dissipation corresponding to the collision operator
    $$
    D_h(g):= -\int Q(g)\log g\,dw \geq 0.
    $$

    \end{enumerate}
\end{lemma}
\begin{proof} We now prove statements (1)-(4) sequentially.  

(1)  For the mass, we simply observe that
    $$
    \rho'(t)=\int L(t)w\cdot\nabla g\,dw
    =-\textnormal{tr}\,L(t)\rho(t).
    $$
    Hence
    $$
    \rho(t)=\rho(0)\exp\left\{-\int_{0}^{t}\textnormal{tr}\,L(s)\,ds\right\}.
    $$
   Since
    $$
\textnormal{tr}\,L(t)=\textnormal{tr}\left(L_0(I+tL_0)^{-1}\right) = \frac{d}{dt}\log \det(I+tL_0),
    $$ 
    we have the second formula for $\rho(t)$: 
    $$
    \rho(0)\exp\left\{-\int_{0}^{t}\textnormal{tr}\,L(s)\,ds\right\} = \frac{\rho(0)}{\det (I+tL_0)}.
    $$
    
  (2)   The evolution of the first moments in \cref{eq:M1} is obtained similarly, by differentiating $\int w_i g(w)\,dw$. 
  
  (3) For the energy, we compute
    $$
    \frac{d}{dt}M_{2}(g) 
    =-\int\textnormal{div}_{w}\left(L(t)w|w|^{2}\right)g
    =-\int\left(2\,L(t)w\cdot w
    +\textnormal{tr}L(t)|w|^{2}\right)g,
    $$
    which gives \cref{eq:second_moment_bound}. 
    
    (4) Finally, the entropy identity \cref{eq:entropy_bound} follows by solving the differential equality
    \begin{align*}
        \frac{d}{dt}\mathcal{H}(g)=\int\left(\log g+1\right)\partial_{t}g & 
        =-D_h(g)+\int\left(\log g+1\right)L(t)w\cdot\nabla g \\
        &=-D_h(g)+\int L(t)w\cdot\nabla\left(g\log g\right) \\
        & =-D_h(g)-\textnormal{tr}\,L(t)\mathcal{H}(g).\qedhere
    \end{align*}
\end{proof}

In the next lemma, we establish estimates for the propagation and generation of moments:

\begin{lemma}\label{propagation_moments}
    Let $g$ be a solution of \cref{eq:drifted_homo_Landau} with $M_s(g_\init)<\infty$, $s\geq 2$, and $\mathcal H(g_\init)<\infty$. Then 
    \begin{equation}\label{eq:prop_moments}
    \sup_{t\in[0,T]}M_{s}(g)\leq C,
    \end{equation}
    where $C$ depends on $s$, $T$, $M_{s}(g_{\init})$, $\mathcal H(g_{\init})$, and $\int_{0}^{T}\left|L(t)\right|\,dt$. 
    
    Moreover, if $\gamma >0$, and $M_\gamma(g_\init)<\infty$, we have that moments of all order $s\geq 2$ become finite at a rate
    \begin{equation}\label{eq:generation_moments}
        M_s(g)\leq Ct^{-\frac{s}{\gamma}},
    \end{equation}
    for small values of $t$, even if $M_s(g_\init)=\infty$.
\end{lemma}
\begin{proof}
    We decompose the equation into its transport and collision parts,
    $$
    \partial_t g=\mathcal Tg+Q(g),
    $$
    and estimate the contribution of each term separately:
    $$
    \frac{d}{dt}M_s(g)=
    \left\langle M_s'(g),\mathcal Tg\right\rangle
    +\left\langle M_s'(g),Q(g)\right\rangle.
    $$
    The transport contribution is easily estimated:
    \begin{align*}
        \left\langle M_{s}'(g),\mathcal{T}g\right\rangle  &
        =\int|w|^{s}\,L(t)w\cdot\nabla_{w}g =-\int\textnormal{div}_{w}\left(L(t)w|w|^{s}\right)g\\
        & =-\int\left(s\,L(t)w\cdot w\, |w|^{s-2}+\textnormal{tr}L(t)|w|^{s}\right)g
        \leq C\left|L(t)\right|M_{s}(g).
    \end{align*}
    The effect of collisions has already been thoroughly studied. For soft potentials
    $-3\le\gamma\le0$, Lemmas 7 and 8 of \cite{CarraDesvHe2017} show that
    $$
    \left\langle M_s'(g),Q(g)\right\rangle
    \le
    \left(D_h(g)+C\right)M_s(g)+C.
    $$
    Combining the previous estimates, we obtain
    $$
    \frac{d}{dt}M_s(g)
    \le \left(D_h(g)+C_T\right)M_s(g)+C.
    $$
    Applying Grönwall's inequality yields
    $$
    M_s(g(t)) \le
    \left(M_s(g_{\init})+Ct\right)
    \exp\left\{
    \int_0^t\left(D_h(g(s))+C\right)\,ds
    \right\},
    $$
    which proves the claim for $-3\le\gamma\le0$ in view of \cref{eq:entropy_bound}.

    For hard potentials, it only remains to estimate the collision term. This was established for
    $0<\gamma\le1$ in \cite{DesvVill2000HardPart1}. We briefly sketch the argument for 
    $\gamma>1$. A straightforward computation shows that
    $$
    \left\langle M_s'(g),Q(g)\right\rangle
    \leq
    s\iint |w-w_*|^\gamma
    \left(
    -|w|^s-|w_*|^s
    +\frac{s}{2}|w|^{s-2}|w_*|^2
    +\frac{s}{2}|w_*|^{s-2}|w|^2
    \right)
    g(w)g(w_*)\,dw\,dw_*.
    $$
    We now bound the right hand side  with a Povzner-type inequality. This inequality was proved in Lemma 2.5 of
    \cite{guo2025kacparticleslandauequation} for $0<\gamma\le 1$ and can be easily extended to all $\gamma >0$, $s>2$ in the form
    \begin{align*}
        & \left(-|w|^{s}-|w_*|^{s}+\frac{s}{2}|w|^{s-2}|w_*|^{2}+\frac{s}{2}|w|^{2}|w_*|^{s-2}\right)|w-w_*|^{\gamma}\\
        \leq\; & -2^{-(\gamma+1)}|w|^{s+\gamma}-2^{-(\gamma+1)}|w_*|^{s+\gamma}+|w|^{s}|w_*|^{\gamma}+|w|^{\gamma}|w_*|^{s}\\
        +\; & C_{\gamma,s}\left(|w|^{s-2+\gamma}|w_*|^{2}+|w|^{2}|w_*|^{s-2+\gamma}\right).
    \end{align*}
    For $s\geq 4$, one can take $C_{\gamma,s}=C_\gamma s^{1+\frac{\gamma}{2}}$.
    This gives
    \begin{equation}\label{eq:moment_bound_big_gamma}
        \left\langle M_{s}'(g),Q(g)\right\rangle \leq -\rho(t)s2^{-\gamma}M_{s+\gamma}(g)+2sM_\gamma(g)M_s(g) + 2sC_{\gamma,s}M_2(g)M_{s-2+\gamma}(g).    
    \end{equation}
    Inequality (\ref{eq:moment_bound_big_gamma})  yields boundedness of $M_s(g)$ for $0<\gamma\le2$, since
    $M_\gamma(g)$ is controlled by $M_2(g)$ and $M_{s-2+\gamma}(g)$ by $M_s(g)$.
    We therefore restrict our attention to the case $\gamma>2$. We first estimate the moment of order $\gamma$ using
    $$
    \left\langle M_\gamma'(g),Q(g)\right\rangle
    \le
    -\rho(t)\gamma2^{-\gamma}M_{2\gamma}(g)
    +2\gamma M_\gamma(g)^2
    +C_\gamma M_2(g)M_{2(\gamma-1)}(g).
    $$
    We estimate $M_{2(\gamma-1)}(g)$ interpolating between $L^1$ and $M_{2\gamma}(g)$, and $M_\gamma(g)$ interpolating between $M_2(g)$ and
    $M_{2\gamma}(g)$:
    \begin{align*}
        C_{\gamma}M_2(g) M_{2(\gamma-1)}(g) & \leq 
        C_{T,\gamma}M_{2\gamma}(g)^{1-\frac{1}{\gamma}} \leq \varepsilon M_{2\gamma}(g)+C_{T,\gamma,\varepsilon};\\
        2\gamma M_\gamma(g)^2 & \leq C_{T,\gamma} M_{2\gamma}(g)^{\frac{\gamma-2}{\gamma-1}} \leq \varepsilon M_{2\gamma}(g) + C_{T,\gamma,\varepsilon}.
    \end{align*}
    Choosing $\varepsilon>0$ sufficiently small yields
    $$
    \frac{d}{dt}M_\gamma(g)\le C_{T,\gamma},
    $$
    which gives a uniform bound on $M_\gamma(g)$. Substituting this estimate into \cref{eq:moment_bound_big_gamma}, we obtain
    \begin{align*}
        \left\langle M_{s}'(g),Q(g)\right\rangle & \leq 
        -\rho(t)s2^{-\gamma}M_{s+\gamma}(g)+C_{T,\gamma}M_s(g) + 2sC_{\gamma,s}M_2(g)M_{s-2+\gamma}(g)\\
        & \leq -\rho(t)s2^{-\gamma}M_{s+\gamma}(g)+C_{T,\gamma}M_s(g) + 2sC_{\gamma,s}M_2(g)\rho(t)^{\frac{2}{s+\gamma}}M_{s+\gamma}(g)^{1-\frac{2}{s+\gamma}}\\
        & \leq -\rho(t)s2^{-(\gamma+1)}M_{s+\gamma}(g)+C_{T,\gamma}M_s(g)+C_{T,\gamma}.
    \end{align*}
    This proves the boundedness of moments. Finally, the generation estimate
    \cref{eq:generation_moments} follows from the inequality
    $$
    M_{s+\gamma}(g)
    \ge
    \rho(t)^{-\frac{\gamma}{s}}
    M_s(g)^{\frac{s+\gamma}{s}},
    $$
    which reduces the problem to solving a differential inequality of the form $\dot y\le -y^{1+\frac{\gamma}{s}}+C$.
\end{proof}

From the previous lemma we get information on the moments and energy for the original function $f$. 

\begin{corollary}\label{cor:moments-f}
    Assume for simplicity that $\int wg_\init(w)\,dw=0$. Then, the bounds (\ref{eq:generation_moments}) and (\ref{eq:prop_moments}) on the moments of $g$ translate  to the following bounds for $f$, solution to \cref{eq:L}: 
    $$
    \int_{\mathbb{R}^3} f(t,x,v)\;dv =  \frac{ \int_{\mathbb{R}^3} f_{\init}(x,v)\;dv}{\det (I+tL_0)},
    $$
    $$
    \sup_{t\in [0,T]}\sup_{x}\int|v-u(t,x)|^sf(t,x,v)\,dv<\infty,
    $$
    where $u(t,x)$ denotes the bulk velocity
    \begin{align*}
    u(t,x):= \frac{1}{\rho} \int_{\mathbb{R}^3} v \; f(t,x,v) \;dv
    =\frac{1}{\rho} \int_{\mathbb{R}^3} (w+ L(t)x) g(w) \;dw
     = L(t)x.
    \end{align*}
     Moreover 
    \begin{align*}
        M_2(f(t,x,\cdot)) = & \; \int_{\mathbb{R}^3} |v|^2f(t,x,v)\;dv = {\overbrace{\int|v-u(t,x)|^2f(t,x,v)\,dv}^{<+\infty}} + |u(t,x)|^2 \rho(t,x) \\
         = & \; C(t) +  |L(t)x|^2 \frac{\rho(0)}{\det(I+tL_0)}. 
    \end{align*}
    \end{corollary}

\vspace{0.5cm}

\begin{remark}
    The bounds of Corollary \ref{cor:moments-f}  show that implosion for the mass and energy of $f$ occurs at the first time $t_*$ such that $\det(I+t_*L_0)=0$. If $L_0 = -I$, one is able to also provide the behavior of the function $f$. This is the content of Theorem \ref{thm:blowup_general}. 
\end{remark}
 
The next lemma shows that the matrix
\begin{align}\label{Diffusion-matrixA[g]}
A[g](w) = \int \Phi(w-z)g(z)\,dz
\end{align}
is uniformly bounded from below.

\begin{lemma}\label{lem:coercivity}
    Let $g$ be a solution of \cref{eq:drifted_homo_Landau} such that $g_\init$ has $M_2(g_\init)<\infty$, $\mathcal H(g_\init)<\infty$ and let $A[g](w)$ be defined as in \eqref{Diffusion-matrixA[g]}.  Then
    $$
    A[g](w)\geq \lambda(t)\langle w\rangle^{\gamma}I,
    $$
    where $I$ denotes the identity, and the constant $\lambda(t)>0$ depends on $\rho(t)$, $M_2(g(t))$, $\mathcal H(g(t))$. 
\end{lemma}
Recall that, by \cref{lem:moments_entropy}, all these quantities remain bounded as long as $\int_0^T|L(t)|\,dt<\infty$.
\begin{proof}
    The bound for soft potentials is a direct consequence of the analysis performed in \cite{Silvestre2015}, where it is shown that if $g$ satisfies estimates of the form
    $$
    M_0 \leq \int g(w)\,dw \leq M_0^{-1},\qquad
    \int |w|^2g(w)\, dw \leq E_0, \qquad
    \int g\log(1+g) \leq H_0,
    $$
    then 
    $$
    A[g](w)\geq \lambda\langle w\rangle^{\gamma}I,
    $$
    where $\lambda$ depends on $M_0$, $E_0$, and $H_0$. An analogous statement for $\gamma\in [0,1]$ can be found in Proposition 4 of \cite{DesvVill2000HardPart1}, which can be generalized to any $\gamma >0$.
\end{proof}

\textbf{Boundedness of Fisher information.}
Boundedness of the Fisher information is the key tool for establishing well-posedness of homoenergetic solutions for very soft potentials. This follows rather directly from the fact that the collision operator dissipates Fisher information. Moreover, the same argument also gives control of the Fisher dissipation functional introduced in \cite{GS24}, namely
\begin{equation}\label{eq:Fisher_dissipation}
    D_{i}(g):=\sum_{k=1}^3\iint\left|\nabla_{w,w_{*}}\left(|w-w_{*}|^{\frac{\gamma}{2}}\tilde{b}_{k}\cdot\nabla_{w,w_{*}}\log G\right)\right|^{2}G\,dw\,dw_{*},\quad G(w,w_{*})=g(w)g(w_{*}),
\end{equation}
where the vector fields $\tilde b_k =\tilde b_k(w-w_*)\in \mathbb R^6$ are defined by
$$
b_{1}(w-w_*):=\begin{pmatrix}0\\
  w^{3}_*-w^{3}\\
  w^{2}-w_*^{2}
  \end{pmatrix},\; b_{2}(w-w_*):=
  \begin{pmatrix}w^{3}-w_*^{3}\\
  0\\
  w_*^{1}-w^{1}
  \end{pmatrix},\; b_{3}(w-w_*):=\begin{pmatrix}w_*^{2}-w^{2}\\
  w^{1}-w_*^{1}\\
  0
  \end{pmatrix},
  $$
  and
  $$
  \tilde{b}_{k}(w-w_*):=\begin{pmatrix}
  b_{k}(w-w_*)\\
  -b_{k}(w-w_*)
  \end{pmatrix},\qquad k=1,2,3.
  $$
   We refer the reader to \cite{Ji2024DisFisher} for lower bounds on $D_i(g)$, which may prove useful for further analysis. Our estimate on the Fisher information can be stated as follows.
\begin{proposition}\label{boundedness_fisher}
    Let $\gamma\in [-3,\sqrt{22}]$, and assume $\det (I+tL_0)>0$ for all $0\leq t\leq T$. Let $g$ be a solution of \cref{eq:drifted_homo_Landau} with $I(g_\init)<\infty$. Then
    
    \begin{align}\label{eq:Fisher}
    e^{-\ell(t)}I\left(g(t)\right)+c_{\gamma}\int_{0}^{t}e^{-\ell(s)}D_{i}\left(g(s)\right)\,ds\leq I(g_{\init}).
    \end{align}
    Here $c_\gamma >0$ for $|\gamma| <\sqrt{22}$ but $c_{\sqrt{22}}=0$.

    The growth rate $\ell(t)$ is: 
    $$\qquad\ell(t):=\int_{0}^{t}\Big(2\left|L(s)\right|-\textnormal{tr}(L(s))\Big)\,ds.
    $$
    
\end{proposition}

\begin{proof}
We decompose the equation into its transport and collision parts,
$$
\partial_{t}g=\mathcal{T}g+Q(g),
$$
so that one can study the effect of each term on the Fisher information:
$$
\frac{d}{dt}I(g)=\left\langle I'(g),\mathcal{T}g\right\rangle +\left\langle I'(g),Q(g)\right\rangle .
$$
In the breakthrough paper \cite{GS24}, followed by \cite{Ji2024DisFisher}, it was shown that the second term is negative, with dissipation given by \cref{eq:Fisher_dissipation}. We thus have
$$
\frac{d}{dt}I(g)\leq-c_\gamma D_i(g)+\left\langle I'(g),\mathcal{T}g\right\rangle ,
$$
and it suffices to bound the effect of the extra transport term on Fisher information. In the following computations, we use an Einstein-like summation convention so that each index should be understood as having a summation in front:
    $$
    \mathcal{T}g=L_{ij}w_{i}\partial_{w_{j}}g.
    $$
    Then we have
    $$
    \left\langle I'(g),\mathcal{T}g\right\rangle 
    =2\int\frac{\partial_{k}g\partial_{k}\left[\mathcal{T}g\right]}{g}-\int\left(\frac{\partial_{k}g}{g}\right)^{2}\mathcal{T}g=:I_{1}+I_{2}.
    $$
    We expand both $I_{1}$, $I_{2}$ and rewrite them in terms of $\sqrt{g}$ using $\nabla g=2\sqrt{g}\nabla\sqrt{g}$.
    \begin{align*}
        I_1  &=  2L_{ij}\int\frac{\partial_{k}g\partial_{k}\left[w_{i}\partial_{j}g\right]}{g}\\
        & =8L_{ij}\int\frac{\partial_{k}\sqrt{g}\partial_{k}\left[w_{i}\sqrt{g}\partial_{j}\sqrt{g}\right]}{\sqrt{g}}\\
        & =8L_{ij}\int\partial_{i}\sqrt{g}\partial_{j}\sqrt{g}+
        8L_{ij}\int\frac{\left(\partial_{k}\sqrt{g}\right)^{2}\partial_{j}\sqrt{g}}{\sqrt{g}}w_{i}+8L_{ij}\int\partial_{k}\sqrt{g}\partial_{kj}\sqrt{g}w_{i}.
    \end{align*}
    Further note that
    \begin{align*}
        8L_{ij}\int\partial_{k}\sqrt{g}\partial_{kj}\sqrt{g}w_{i} &
        = 4L_{ij}\int\partial_{j}\left[\left(\partial_{k}\sqrt{g}\right)^{2}\right]w_{i}\\
        & =-4L_{ii}\int\left(\partial_{k}\sqrt{g}\right)^{2}=-\,\textnormal{tr}(L)I(g).
    \end{align*}
    Therefore,
    $$
    I_{1}=8\int L\nabla\sqrt{g}\cdot\nabla\sqrt{g} -\,\textnormal{tr}(L)I(g)+8L_{ij}\int\frac{\left(\partial_{k}\sqrt{g}\right)^{2}\partial_{j}\sqrt{g}}{\sqrt{g}}w_{i}.
    $$
    As for the second term $I_2$,
    $$
    I_{2}=-L_{ij}\int\left(\frac{\partial_{k}g}{g}\right)^{2}w_{i}\partial_{w_{j}}g=
    -8L_{ij}\int\frac{\left(\partial_{k}\sqrt{g}\right)^{2}\partial_{j}\sqrt{g}}{\sqrt{g}}w_{i},
    $$
    which cancels with the last term of $I_{1}$. All in all, we obtain
    $$
    \left\langle I'(g),\mathcal{T}g\right\rangle =8\int L\nabla\sqrt{g}\cdot\nabla\sqrt{g} -\,\textnormal{tr}(L)I(g)\leq\Big(2\left|L\right|-\textnormal{tr}\,L\Big)I(g).
    $$
    This shows the differential inequality
    \begin{equation}\label{eq:diff_ineq_Fisher}
        \frac{d}{dt} I (g) \leq -c_\gamma D_i(g) + \Big(2\left|L(t)\right|-\textnormal{tr}(L(t))\Big)I(g),
    \end{equation}
    which implies the stated estimate.
\end{proof}

    
    The qualitative behavior of the Fisher information of $g$ will strongly depend on the deformation matrix \(L(t)\). Consider the differential inequality \cref{eq:diff_ineq_Fisher}. If \(L_0 = \textnormal{diag}(\lambda_1,\lambda_2,\lambda_3)\) with $\lambda_i>0$, then as \(t\to\infty\), we have
    $$
    L(t) = \frac{1}{t}I+O\left(\frac{1}{t^2}\right),
    $$
    which gives
    $$
    \frac{d}{dt}I(g)
    \leq
    \Bigl(2|L(t)|-\operatorname{tr}L(t)\Bigr)I(g)
    =
    \left(\frac{2}{t}-\frac{3}{t}+O(t^{-2})\right)I(g).
    $$
    Thus, after a transient period, the Fisher information becomes monotone decreasing. By contrast, if \(L_0\) corresponds to a simple shear, see \cref{eq:shear_mat}, there is no apparent reason to expect such monotonicity.
    
    A similar phenomenon occurs for the second moment \(M_2(g)\). In the case of homogeneous dilation it is monotone decreasing, whereas this property generally fails for an arbitrary deformation matrix \(L(t)\). This suggests that a more refined qualitative analysis of homoenergetic solutions should be carried out on a case-by-case basis.

\begin{remark}
   Estimate (\ref{eq:Fisher}) yields the following bound for any  $f$ solution to \cref{eq:L}: 
    $$
    \sup_{x\in\mathbb R^3}    \int \frac{|\nabla_v f|^2}{f} \,dv = I(g(t)),\qquad 
    \sup_{x\in\mathbb R^3} \int \frac{|\nabla_x f|^2}{f} \,dv\leq |L(t)|^2I(g(t)).
    $$
    \end{remark}

\subsection{Global-in-time existence of smooth solution}
The main goal of this section is to prove the following $L^\infty$ estimate, since higher regularity theory after that point is quite standard.
\begin{theorem}\label{thm:boundedness}
    Assume $\det(I+tL_0)>0$ for all $0\leq t\leq T$. Let $g$ be a solution of \cref{eq:drifted_homo_Landau} with $g_\init\in L^\infty\cap L^1$, $M_2(g_\init)<\infty$. Further assume
    \begin{enumerate}
        \item In the soft potential case $-3\leq \gamma <0$, that $I(g_\init)<\infty$ and $M_{s}(g_\init)<\infty$ for some $s > \max\left\{ \frac{2|\gamma|}{4+\gamma},|\gamma|,2\right\}$. Then
        \begin{equation*}
            \sup_{0\leq t\leq T}\Vert g(t)\Vert_{L^\infty} \leq C,
        \end{equation*}
        where $C$ depends on $\gamma$, $T$, $\int_0^T|L(t)|\,dt$, $s$, $\rho(0)$, $M_{s}(g_\init)$, $I(g_\init)$, and $\Vert g_\init\Vert_{L^\infty}$.
        \item In the hard potential case $\gamma \geq 0$, that $M_{\frac{3\gamma}{2}}(g_\init)<\infty$. Then
        \begin{equation*}
            \sup_{0\leq t\leq T}\Vert g(t)\Vert_{L^\infty} \leq C,
        \end{equation*}
        where $C$ depends on $\gamma$, $T$, $\int_0^T|L(t)|\,dt$, $\rho(0)$, $M_2(g_\init)$, $M_{\frac{3\gamma}{2}}(g_\init)$, and $\Vert g_\init\Vert_{L^\infty}$.
    \end{enumerate}
\end{theorem}
Note that $\mathcal H(g_\init)$ is automatically finite under the previous assumptions.
The proof is based on a Moser iteration which we will perform for the whole range $\gamma\geq-3$. We require several lemmas, beginning with the following $\varepsilon$-Poincaré inequality:
\begin{lemma}\label{eps_poincare}
    Let $-3\leq \gamma <0$ and denote
    \begin{equation}\label{eq:cgamma}
    c_{\gamma}[g](w):=\begin{cases}
    2(3+\gamma)\big(g*|\cdot|^{\gamma}\big)(w), 
    & -3<\gamma<0\\
    8\pi g, & \gamma=-3
    \end{cases}.
    \end{equation}
    Assume that $g$ has finite mass, $I(g)<\infty$, and $M_{\mathrm{m}}(g)<\infty$  with $\mathrm{m}>\max\left\{ \frac{2|\gamma|}{4+\gamma},|\gamma|\right\}$ . Then, for any $0<\varepsilon<1$ we have
    $$
    \int g^{p}c_{\gamma}[g]\leq
    \varepsilon\int|\nabla g^{\frac{p}{2}}|^{2}\langle w\rangle^{\gamma}+C\varepsilon^{-M}\int g^{p}\langle w\rangle^{\gamma},
    $$
    where $M$ depends on $\gamma$, and $C$ further depends on $\mathrm{m}$, the mass of $g$, $M_\mathrm{m}(g)$, $I(g)$.
\end{lemma}
\begin{proof}
    The starting point is the $\varepsilon$-Poincaré inequality of Proposition 1.7 of \cite{ABDL24}, which states that
    \begin{equation}\label{eq:alonso_eps_Poinc}
        \int\phi^{2}c_\gamma[g]\leq\varepsilon\int\left|\nabla\left(\langle w\rangle^{\frac{\gamma}{2}}\phi\right)\right|^{2}+C_{0}\left(\left\Vert g\right\Vert _{L^{1}}+\varepsilon^{-\frac{s}{1-s}}\left\Vert \langle\cdot\rangle^{-\gamma}g\right\Vert _{L^{q}}^{\frac{1}{1-s}}\right)\int\phi^{2}\langle w\rangle^{\gamma},
    \end{equation}
    where $\max\left\{ \frac{3}{5+\gamma},1\right\} <q<\frac{3}{3+\gamma}$, and $s=\frac{3-q(3+\gamma)}{2q}\in (0,1)$. Now we bring the factor $\langle w\rangle^{\frac{\gamma}{2}}$ out of the gradient, which will make the inequality more convenient to use in the sequel
    \begin{align*}
        \int\left|\nabla\left(\langle w\rangle^{\frac{\gamma}{2}}\phi\right)\right|^{2}& 
        =\int\left|\langle w\rangle^{\frac{\gamma}{2}}\nabla\phi+\phi\nabla\langle w\rangle^{\frac{\gamma}{2}}\right|^{2}\\
        & =\int\left|\nabla\phi\right|^{2}\langle w\rangle^{\gamma}+2\int\left(\nabla\phi\cdot\nabla\langle w\rangle^{\frac{\gamma}{2}}\right)\langle w\rangle^{\frac{\gamma}{2}}\phi+\int\left|\nabla\langle w\rangle^{\frac{\gamma}{2}}\right|^{2}\phi^{2}\\
        & =\int\left|\nabla\phi\right|^{2}\langle w\rangle^{\gamma}+\int\left(\left|\nabla\langle w\rangle^{\frac{\gamma}{2}}\right|^{2}-\frac{1}{2}\Delta\langle w\rangle^{\gamma}\right)\phi^{2}\\
        & \leq\int\left|\nabla\phi\right|^{2}\langle w\rangle^{\gamma}+15\int\phi^{2}\langle w\rangle^{\gamma}.
    \end{align*}
    Plugging this into \cref{eq:alonso_eps_Poinc} and setting $\phi=g^{\frac{p}{2}}$ yields
    $$
    \int g^{p}c_\gamma[g]\leq
    \varepsilon\int\left|\nabla\phi\right|^{2}\langle w\rangle^{\gamma}+
    C_{0}\left(\left\Vert g\right\Vert _{L^{1}}+\varepsilon^{-\frac{s}{1-s}}\left\Vert \langle\cdot\rangle^{-\gamma}g\right\Vert _{L^{q}}^{\frac{1}{1-s}}+15\right)
    \int g^{p}\langle w\rangle^{\gamma}.
    $$
    We already have control of the mass of $g$, so it is only left to show that $\left\Vert \langle\cdot\rangle^{-\gamma}g\right\Vert _{L^{q}}$ is bounded. This can be done by considering
    \begin{equation}\label{eq:constraint_q}
        \max\left\{ \frac{3}{5+\gamma},1\right\} <q<\min\left\{ 3,\frac{3}{3+\gamma}\right\} 
    \end{equation}
    and interpolating $\left\Vert \langle\cdot\rangle^{-\gamma}g\right\Vert _{L^{q}}$ between $L^3$ and $L^1_\mathrm{m}$, using that the $L^3$ norm of $g$ is bounded by its Fisher information thanks to the Sobolev inequality
    $$\|g\|_{L^3}=\|\sqrt g\|_{L^6}^2\le C_{\rm Sob}^2\|\nabla\sqrt g\|_{L^2}^2=\tfrac{C_{\rm Sob}^2}{4}\,I(g).$$
    A simple computation shows that for this choice to be possible, one only needs to have bounded moments of order $\mathrm{m}>\max\left\{ \frac{2|\gamma|}{4+\gamma},|\gamma|\right\}$, which finishes the proof.
\end{proof}

Next, we perform a simple $L^p$ energy estimate.

\begin{lemma}\label{lem:Lp_energy_estimate}
    Assume the setting of \cref{thm:boundedness} and let $p\geq 2$, $0\leq t\leq T$.
    \begin{enumerate}[label=\roman*)]
        \item If $\gamma\geq0$,
        $$
        \frac{d}{dt}\int g^{p}+\lambda\int\langle w\rangle^{\gamma}\left|\nabla g^{\frac{p}{2}}\right|^{2}\leq Cp\int\langle w\rangle^{\gamma}g^{p}.
        $$
        \item If $\gamma<0$,
        $$
        \frac{d}{dt}\int g^{p}+\lambda\int\langle w\rangle^{\gamma}\left|\nabla g^{\frac{p}{2}}\right|^{2}\leq C\int g^{p}+Cp^{M}\int\langle w\rangle^{\gamma}g^{p}.
        $$
    \end{enumerate}
    Here $C$ has the dependence described in \cref{thm:boundedness}, and $M$ depends only on $\gamma$. Crucially, $C$ does not depend on $p$.
\end{lemma}
\begin{proof}
    We start with a direct computation $$
    \frac{d}{dt}\int g^{p}=p\int g^{p-1}\left( L(t)w\cdot\nabla g\right) +p\int g^{p-1}Q(g).
    $$
    The contribution of the drift term is
    $$
    p\int g^{p-1}\left( L(t)w\cdot\nabla g\right) =\int\left( L(t)w\cdot\nabla g^{p}\right) =-\textnormal{tr}\,L(t)\int g^{p}.
    $$
    To estimate the contribution of the collision operator, we write it as
    $$
    Q(g)=\textnormal{div}\left(A[g]\nabla g-b[g]g\right),\qquad 
    b[g](w):=-2\int\left|w-w_*\right|^{\gamma}(w-w_*)g(w_*)\,dw_*.
    $$
    Then
    \begin{align*}
        p\int g^{p-1}Q(g) & 
        =-p(p-1)\int\left( A[g]\nabla g\cdot\nabla g\right) g^{p-2}+
        p(p-1)\int b[g]\cdot\nabla g\,g^{p-1}\\
        & =-\frac{4(p-1)}{p}\int\left( A[g]\nabla g^{\frac{p}{2}}\cdot\nabla g^{\frac{p}{2}}\right) -(p-1)\int\textnormal{div}\,b[g]\,g^{p}.
    \end{align*}
    Now we use that
    $$
    -\textnormal{div}\,b[g] = c_\gamma[g] =
    \begin{cases}
        2(\gamma+3)|w|^{\gamma}*g, & \gamma>-3\\
        8\pi g, & \gamma = -3.
    \end{cases},
    $$
    where $c_\gamma[g]$ was defined in \cref{eps_poincare} for $\gamma<0$, although the previous relation holds for all $\gamma\geq -3$. Using the ellipticity estimate from \cref{lem:coercivity} and putting everything together, we obtain
    $$
    \frac{d}{dt}\int g^{p}+2\lambda\int\langle w\rangle^{\gamma}\left|\nabla g^{\frac{p}{2}}\right|^{2}\leq C\int g^{p}+p\int c_{\gamma}[g]g^{p}.
    $$
    \textit{Case of hard potentials $\gamma\geq0$.} We are done, since
    $$
    c_{\gamma}[g]=\int|w-w_{*}|^{\gamma}g(w_{*})\leq C\int\left(|w|^{\gamma}+|w_{*}|^{\gamma}\right)g(w_{*})\leq C\langle w\rangle^{\gamma}.
    $$
    \textit{Case of soft potentials.} Using the $\varepsilon$-Poincaré inequality from \cref{eps_poincare} gives
    $$
    \frac{d}{dt}\int g^{p}+2\lambda\int\langle w\rangle^{\gamma}\left|\nabla g^{\frac{p}{2}}\right|^{2}\leq 
    C\int g^{p}+p\varepsilon\int|\nabla g^{\frac{p}{2}}|^{2}\langle w\rangle^{\gamma}+Cp\varepsilon^{-M}\int g^{p}\langle w\rangle^{\gamma}.
    $$
    Then, choosing $\varepsilon=\lambda/p$ finishes the proof.
\end{proof}

As an immediate consequence, we have control of $L^p$ norms for any finite $p$:

\begin{corollary}\label{lem:boundedness_Lp}
    Assume the setting of \cref{thm:boundedness} and let $p\geq 2$.
    \begin{enumerate}[label=\roman*)]
        \item For $\gamma \geq 0$,
        $$
        \sup_{t\in[0,T]}\left\Vert g\right\Vert _{L^{p}} +\left(\int_{0}^{T}\int\langle w\rangle^{\gamma}\left|\nabla g^{\frac{p}{2}}\right|^{2}\right)^{\frac{1}{p}}
        \leq C\left(\left\Vert g_{\init}\right\Vert _{L^{p}}+p^{\frac{3}{2}}M_{\frac{3}{2}\gamma}(g_{\init})+1\right).
        $$
        \item For $\gamma <0$,
        $$
        \sup_{t\in[0,T]}\left\Vert g\right\Vert _{L^{p}}+
        \left(\int_{0}^{T}\int\langle w\rangle^{\gamma}\left|\nabla g^{\frac{p}{2}}\right|^{2}\right)^{\frac{1}{p}}\leq
        C^{p^{M}}\left\Vert g_{\init}\right\Vert _{L^{p}}.
        $$
    \end{enumerate}
    The constant $C$ has the same dependence as in \cref{thm:boundedness}, and $M$ depends only on $\gamma$.
\end{corollary}
\begin{proof}
    The estimate is immediate for soft potentials. For hard potentials, we use Hölder inequality with exponents $\frac{3p-1}{p-1}$, $\frac{3p-1}{2p}$ with respect to the measure $g$
    $$
    pC\int\langle w\rangle^{\gamma}g^{p}\leq Cp\left(\int|w|^{\gamma\frac{3p-1}{2p}}g\right)^{\frac{2p}{3p-1}}\left(\int g^{3p}\right)^{\frac{p-1}{3p-1}}.
    $$
    Next, we note that $\frac{p-1}{3p-1}<\frac{1}{3}$, so we can use Young's inequality with exponents $\frac{3p-1}{3p-3}$, $\frac{3p-1}{2}$ together with Sobolev's inequality to obtain
    \begin{align*}
        pC\int|w|^{\gamma}g^{p}&
        \le\varepsilon\left(\int g^{3p}\right)^{\frac{1}{3}}+
        C^{\frac{3p-1}{2}}\varepsilon^{-(p-1)}p^{\frac{3}{2}p}M_{\gamma\frac{3p-1}{2p}}(g)^{p}\\
        &\leq\varepsilon C_{\textnormal{Sob}}\int\langle w\rangle^{\gamma}\left|\nabla g^{\frac{p}{2}}\right|^{2}+
        C^{\frac{3p-1}{2}}\varepsilon^{-(p-1)}p^{\frac{3}{2}p}M_{\gamma\frac{3p-1}{2p}}(g)^{p}.
    \end{align*}
    Now the proof is finished choosing $\varepsilon=\frac{\lambda}{2C_{\textnormal{Sob}}}$ and using \cref{lem:Lp_energy_estimate}.
\end{proof}

Now we are ready to prove the main result of the section.

\begin{proof}[Proof of \cref{thm:boundedness}]
    Our main tool is the following inequality, which follows by applying Gronwall's inequality on \cref{lem:Lp_energy_estimate}:
    \begin{equation}\label{eq:Lp_Lpgamma}
        \sup_{t\in[0,T]}\int g^{p}+\int_{0}^{T}\int\langle w\rangle^{\gamma}\left|\nabla g^{\frac{p}{2}}\right|^{2}\leq C\left(\int g_{\init}^{p}+p^{M}\int_{0}^{T}\int\langle w\rangle^{\gamma}g^{p}\right).
    \end{equation}
    The constant $C$ is uniform in $p$, and has the dependence described in \cref{thm:boundedness}. Moreover, we stress that the previous inequality is valid for all $\gamma \geq -3$. In fact, the Moser iteration takes the same form for hard and soft potentials. We will also use the Sobolev inequality
    \begin{equation}\label{eq:Sobolev}
        \left(\int\langle w\rangle^{3\gamma}g^{3p}\right)^{\frac{1}{3}}
        \leq C\int\left|\nabla\left(\langle w\rangle^{\frac \gamma 2}g^{\frac{p}{2}}\right)\right|^{2}
        \leq C\int\langle w\rangle^{\gamma}\left|\nabla g^{\frac{p}{2}}\right|^{2}+C\int\langle w\rangle^{\gamma}g^{p}.
    \end{equation}
    Together, these inequalities provide control of $\sup_{t\in[0,T]}\int g^{p}$ and $\int_0^T\int\langle w\rangle^{3\gamma}g^{3p}$. Then, $L_{t,w}^{\frac{5}{3}p}$ appears as the natural space to do the Moser iteration. To be precise, using Hölder's inequality with exponents $(3,\frac{3}{2})$ gives
    \begin{align*}
        \int_{0}^{T}\int\langle w\rangle^{\gamma}g^{\frac{5}{3}p}&
        =\int_{0}^{T}\int\langle w\rangle^{\gamma}g^{p}\cdot g^{\frac{2}{3}p}\\
        &\leq\int_{0}^{T}\left(\int\langle w\rangle^{3\gamma}g^{3p}\right)^{\frac{1}{3}}\left(\int g^{p}\right)^{\frac{2}{3}}\,dt\\
        &\leq\sup_{t\in[0,T]}\left(\int g^{p}\right)^{\frac{2}{3}}\int_{0}^{T}\left(\int\langle w\rangle^{3\gamma}g^{3p}\right)^{\frac{1}{3}}.
    \end{align*}
    Now, using \cref{eq:Lp_Lpgamma,eq:Sobolev} yields
    \begin{align*}
        \int_{0}^{T}\int\langle w\rangle^{\gamma}g^{\frac{5}{3}p}&
        \leq C\sup_{t\in[0,T]}\left(\int g^{p}\right)^{\frac{2}{3}}\left(\int_{0}^{T}\int\langle w\rangle^{\gamma}\left|\nabla g^{\frac{p}{2}}\right|^{2}
        +\int_{0}^{T}\int\langle w\rangle^{\gamma}g^{p}\right)\\
        &\leq C\left(\int g_{\init}^{p}+p^M\int_{0}^{T}\int\langle w\rangle^{\gamma}g^{p}\right)^{\frac{5}{3}}.
    \end{align*}
    Raising to $\frac{3}{5p}$, using $\left(\int g_\init^p\right)^{\frac1p}\leq \max\{1,\rho(0)\}^{\frac{1}{p}}\left\Vert g_{\init}\right\Vert _{L^{\infty}}^{1-\frac{1}{p}}$, and absorbing the factor $\max\{1,\rho(0)\}^{\frac1p}$ into $C^{\frac1p}$, we obtain
    $$
    \left(\int_{0}^{T}\int\langle w\rangle^{\gamma}g^{\frac{5}{3}p}\right)^{\frac{3}{5p}}\leq C^{\frac{1}{p}}p^{\frac{M}{p}}
    \max\left\{1, \left\Vert g_{\init}\right\Vert _{L^{\infty}},\left(\int_{0}^{T}\int\langle w\rangle^{\gamma}g^{p}\right)^{\frac{1}{p}}\right\} .
    $$
    Here, if we define
    $$
    p_{k}=2\left(\frac{5}{3}\right)^{k},\qquad 
    E_{k}=\left(\int_{0}^{T}\int\langle w\rangle^{\gamma}g^{p_{k}}\right)^{\frac{1}{p_{k}}},
    $$
    we obtain the recursive bound
    $$
    E_{k+1}\leq 
    C^{\frac{1}{p_{k}}}p_{k}^{\frac{M}{p_{k}}}\max\left\{1, \left\Vert g_{\init}\right\Vert _{L^{\infty}},E_{k}\right\} .
    $$
    This implies
    $$
    \left\Vert g\right\Vert _{L^{\infty}\left([0,T]\times\mathbb{R}^{3}\right)}=E_{\infty}\leq
    \prod_{k=0}^{+\infty}C^{\frac{1}{p_{k}}}p_{k}^{\frac{M}{p_{k}}}\cdot\max\left\{1, \left\Vert g_{\init}\right\Vert _{L^{\infty}},E_{0}\right\},
    $$
    which reduces the problem to simply estimating
    $$
    E_0^2 = \int_{0}^{T}\int\langle w\rangle^{\gamma}g^{2}.
    $$
    For soft potentials, this is a consequence of \cref{lem:boundedness_Lp}, since $\langle w\rangle^\gamma\leq 1$. For hard potentials, H\"older's inequality gives
    $$
    \int\langle w\rangle^{\gamma}g^{2}\leq\left(\int g^{4}\right)^{\frac{1}{3}}\left(\int\langle w\rangle^{\frac{3\gamma}{2}}g\right)^{\frac{2}{3}},
    $$
    and the right hand side is bounded uniformly in $t\in[0,T]$ by \cref{lem:boundedness_Lp} with $p=4$ and \cref{propagation_moments}.
\end{proof}

\textit{A comment on uniqueness.}
Uniqueness of solutions to \cref{eq:drifted_homo_Landau} can be obtained by adapting
the techniques developed for the spatially homogeneous equation. In the case of soft
potentials $-3 \leq \gamma \leq 0$, for instance, the arguments of
\cite{fournier2009well_posedness_soft_potentials,fournier2010uniqueness_Coulomb,Caja2026uniquenessFuzzy}
apply with minor modifications and yield uniqueness of weak solutions to
\cref{eq:drifted_homo_Landau} within the class $L^1_t L^{p_\gamma}_{w}(\mathbb{R}^3)$,
where $p_\gamma = \tfrac{3}{3+\gamma} \in [1,\infty]$. 

For hard potentials we refer to
\cite{FournierHeydecker2021}, where uniqueness of weak solutions is established under moment assumptions. By contrast, uniqueness for the inhomogeneous equation \cref{eq:L} is considerably more delicate, and we are not aware of any result that applies in our setting beyond the class of bounded smooth solutions.

\section{Hydrodynamic blow-up for singular deformation matrices}\label{sec:density_blow_up}
We now consider the case in which
\[
t_*:=\inf\{t>0:\det(I+tL_0)=0\}<\infty.
\]
This occurs precisely when \(L_0\) has a real negative eigenvalue; if \(\lambda<0\) is the eigenvalue corresponding to the first singular time, then \(t_*=-\lambda^{-1}\). By \cref{lem:moments_entropy}, the mass of \(g\) and density of \(f\) are given by
\[
\rho(t)=\frac{\rho(0)}{\det(I+tL_0)} \gtrsim |t-t_*|^{-n_\lambda},
\]
where \(n_\lambda\) is the multiplicity of the eigenvalue \(\lambda\). Although the argument is elementary, the result is important: to the best of our knowledge, it provides the first example of finite-time singularity formation for the inhomogeneous Landau equation for $\gamma\leq \sqrt 3$.

A natural next step is to study the behavior of the solution as \(t\to t_*\). In particular, we ask:
\begin{center}
\itshape
Is the singularity confined to the blow-up of the hydrodynamic quantities of \(f\), or is it also accompanied by a loss of regularity?
\end{center}

As a first step toward this question, we focus on \(L_0=-I\), leaving the more general case for future work. This corresponds to initial data of the form $f_\init(x,v)=g_\init(x+v)$. Computations in \cref{sec:derivation} (with $L_0=-I$, and $\xi(t,x)=L(t)x=-\frac{x}{1-t}$) show that the solution of (\ref{eq:L}) takes the form
$$
f(t,x,v)=g\left(t,v+\frac{x}{1-t}    \right). 
$$
Motivated by this observation we perform a self-similar ansatz
$$
f(t,x,v)=F\left(s(t),\frac{x}{(1-t)^{c_x}},\frac{v}{(1-t)^{c_v}}\right),
\qquad
X:=\frac{x}{(1-t)^{c_x}},\quad V:=\frac{v}{(1-t)^{c_v}},
$$
where \(s(t)\) is a time reparametrization to be chosen below.
A direct computation gives
\begin{align*}
    \partial_tf
    &=(1-t)^{-1}\Big[(1-t)s'(t)\partial_sF
    +c_xX\cdot\nabla_XF+c_vV\cdot\nabla_VF\Big],\\
    v\cdot \nabla_xf
    &=(1-t)^{c_v-c_x}V\cdot \nabla_XF,\\
    Q(f)
    &=(1-t)^{c_v(3+\gamma)}Q(F).
\end{align*}
Imposing 
$$
c_v-c_x=-1
$$
matches the scaling of the transport terms, leading to
$$
(1-t)s'(t)\partial_sF
+c_x X\cdot \nabla_XF
+V\cdot(\nabla_X+c_v\nabla_V)F
=(1-t)^{c_v(3+\gamma)+1}Q(F).
$$
Now there are two natural choices for the parameters $c_x,c_v$. The choice $c_v=0,c_x=1$ corresponds to the computations in \cref{sec:derivation} and leads to the homogeneous equation with an extra hyperbolic term. One can also choose 
$$
c_v = -1,\qquad c_x = 0;
$$
this choice allows us to connect (\ref{eq:L}) with the homogeneous equation. Indeed, 
\begin{align*}
    f(t,x,v) = F\left(s(t),x,(1-t)v\right)= F(s(t), X,V),
\end{align*}
and $F$ solves 
\begin{equation}\label{eq:self_similar2}
    (1-t)s'(t)\partial_sF
+V\cdot(\nabla_X-\nabla_V)F
=(1-t)^{-2-\gamma}Q(F).
\end{equation}
We look for solutions of the form 
$$
F(s,X,V)=F(s,X+V),
$$
and \cref{eq:self_similar2} becomes
$$
s'(t)\partial_s F = (1-t)^{-3-\gamma}Q(F).
$$
Choosing the reparameterization
$$
s(t) = \int_0^t(1-\tau)^{-3-\gamma}\,d\tau,
$$
we obtain that
\begin{equation}\label{eq:final_self_similar}
    f(t,x,v) = F(s(t),W),\qquad W = x+(1-t)v,
\end{equation}
where now $F$ solves the homogeneous Landau equation
$$
\begin{cases}
    \partial_sF=Q(F),\\
    F(0,W)=g_\init(W).
\end{cases}
$$
With this setup, studying the behavior of $f$ as $t\to 1^{-}$ becomes straightforward.
\begin{proof}[Proof of \cref{thm:blowup_general}]
    The fact that $f\in \mathcal C^\infty\left([0,1)\times \mathbb R^6\right)\cap L^\infty\left([0,1)\times \mathbb R^6\right)$ is a direct consequence of \cref{eq:final_self_similar} and the fact that $F$ solves the homogeneous Landau equation, which we know is well posed. Moreover, \cref{implos_rho,implos_theta,implos_u} now follow by a direct computation:
    \begin{align*}
        \rho(t,x)&=\int F\left(s(t),x+(1-t)v\right)\,dv\\
        &=(1-t)^{-3}\int F\left(s(t),W\right)\,dW=(1-t)^{-3}\int g_{\init}(W)\,dW.
    \end{align*}
    Similarly, for the bulk velocity we have
    \begin{align*}
        u(t,x)&=\frac{1}{\rho(t,x)}\int vF\left(s(t),x+(1-t)v\right)\,dv\\
        &=\frac{(1-t)^{-4}}{\rho(t,x)}\int(W-x)F\left(s(t),W\right)\,dW\\
        &=-\frac{x}{(1-t)^{4}}\frac{\rho(0)}{\rho(t)}=-\frac{x}{1-t},
    \end{align*}
    where we used that $\int WF(W)\,dW=0$. Finally, for the temperature we have
    \begin{align*}
        \theta(t,x)&=\frac{1}{3\rho(t)}\int|v-u(t,x)|^{2}F\left(s(t),x+(1-t)v\right)\,dv\\
        &=\frac{1}{3\rho(t)(1-t)^{2}}\int|(1-t)v+x|^{2}F\left(s(t),x+(1-t)v\right)\,dv\\
        &=\frac{1}{3\rho(t)(1-t)^{5}}\int|W|^{2}F(s(t),W)\,dW\\
        &=\frac{1}{3}(1-t)^{-2}\frac{\int|W|^{2}g_{\init}(W)\,dw}{\int g_{\init}(W)\,dw},
    \end{align*}
    which finishes the proof.
\end{proof}
We recall that, in \cref{thm:blowup_general}, the solution takes the form \cref{eq:final_self_similar}. Written more explicitly:
\begin{equation}\label{eq:self_similar_final}
    f(t,x,v)=F\big(s(t),x+(1-t)v\big),\qquad
s(t)=\int_0^t(1-\tau)^{-\gamma-3}\,d\tau.
\end{equation}
Note that, if we fix $x$, as $t\to1^{-}$, the function $f(t,x,\cdot)$ is approaching the constant value $F\big(s(1^{-}),x\big)$. Therefore, the imploding behavior we observe in $(\rho,u,\theta)$ is caused by tail fattening. On the other hand, the shape of the profile $F(s(t),\cdot)$ depends on the parameter~$\gamma$:
\begin{itemize}
    \item If $\gamma <-2$, then $s(1^{-})<\infty$, so $F\big(s(1^{-}),\cdot\big)$ corresponds to an intermediate profile between $g_\init$ and a global Maxwellian. For example, if $\gamma = -3$ then $s(t)=t$ and~$s(1^{-}) = 1$.
    \item If $\gamma \geq -2$, then $s(1^{-})=\infty$, so $F\big(s(1^{-}),\cdot\big)$ is in fact a global Maxwellian. This shows that for $t\sim 1$, the solution \cref{eq:self_similar_final} behaves approximately as
    \begin{align*}
        f(t,x,v) & \sim \frac{\rho(0)}{(2\pi\theta(0))^{\frac{3}{2}}}
        \exp\left\{-\frac{|x+(1-t)v|^2}{2\theta(0)}\right\}\\
        & =\frac{\rho(t)}{(2\pi\theta(t))^{\frac{3}{2}}}\exp\left\{-\frac{|u(t,x)-v|^2}{2\theta(t)}\right\}.
    \end{align*}
\end{itemize}
The distinction of the asymptotic profile between very soft and $\gamma \geq-2$ is quite interesting. For nonsingular kernels of $Q$ ($\gamma >-2$), this type of blow-up profile relates to a Maxwellian function. However, if $\gamma<-2$, the blow-up profile is given by the solution to the homogeneous equation at some intermediate time. After the completion of this work, we discovered that this distinction of the asymptotic profile based on the value of $\gamma$, was also identified by Nikolskii in \cite{Nikolskii1963} for the Boltzmann equation.

\appendix

\section{Classification of long time asymptotics for deformation matrices}\label{sec:asymptotics}

James, Nota, and Vel\'azquez classified in \cite{james2019self} all possible asymptotic behaviors as $t\to +\infty$ of the deformation matrix 
$$
L(t):= L_0(I+tL_0)^{-1}
$$ 
for any matrix $L_0$ such that 
$$
\det(I+tL_0) > 0, \qquad \textnormal{for all }t\geq0.
$$
 Here we provide the complete list:
\begin{itemize}
    \item Homogeneous dilatation 
    $$
    L(t)=\frac{1}{t}I+O\left(\frac{1}{t^{2}}\right),\qquad t\to\infty.
    $$
    \item Cylindrical dilatation and shear
    $$
    L(t)=\frac{1}{t}\begin{pmatrix}
    1 & 0 & K\\
    0 & 1 & 0\\
    0 & 0 & 0
    \end{pmatrix}+O\left(\frac{1}{t^{2}}\right),\qquad t\to\infty,
    $$
    when $K\neq0$ or simple cylindrical dilatation when $K=0$.
    \item Planar shear
    $$
    L(t)=\frac{1}{t}\begin{pmatrix}0 & 0 & 0\\
    0 & 0 & K\\
    0 & 0 & 1
    \end{pmatrix}+O\left(\frac{1}{t^{2}}\right),\qquad t\to\infty.
    $$
    \item Simple shear
    $$
    L(t)=\begin{pmatrix}0 & K & 0\\
    0 & 0 & 0\\
    0 & 0 & 0
    \end{pmatrix} = L_0,\qquad K\neq0.
    $$
    \item Simple shear with decaying planar dilatations/shear:
    $$
    L(t)=\begin{pmatrix}0 & K_{2} & 0\\
    0 & 0 & 0\\
    0 & 0 & 0
    \end{pmatrix}+\frac{1}{t}\begin{pmatrix}0 & K_{1}K_{3} & K_{1}\\
    0 & 0 & 0\\
    0 & K_{3} & 1
    \end{pmatrix}+O\left(\frac{1}{t^{2}}\right),\qquad K_{2}\neq0,\qquad t\to\infty.
    $$
    \item Combined orthogonal shear:
    $$
    L(t)=\begin{pmatrix}0 & K_{3} & K_{2}-tK_{1}K_{3}\\
    0 & 0 & K_{1}\\
    0 & 0 & 0
    \end{pmatrix},\qquad K_{1}K_{3}\neq0,\qquad t\to\infty.
    $$
\end{itemize}

\medskip

\textbf{A comment on the expected asymptotics.}

The long-time behavior of solutions to \eqref{eq:drifted_homo_Landau} can be understood comparing the order of magnitude of the hyperbolic term $L(t)w\cdot\nabla_w g$ against that of $Q(g)$, both expressed through the order-of-magnitude scales $\rho(t)$ (density), $[w](t)$ (typical velocity scale), and $[L](t)$ (scale of $L(t)$).

For the Landau equation, the collision operator $Q$ is quadratic and rescales as:
\begin{equation}
\rho(t) [w]^{\gamma} [g], \label{eq:colltermscaling}
\end{equation}
where $[g]$ denotes the order of magnitude of $g$, and $\gamma$ parameterizes the collision kernel $\Phi$ (assumed to be homogeneous of degree $\gamma+2$). With convention, $\gamma$ directly corresponds to the homogeneity parameter in the Boltzmann collision kernel.  
The deformation term $L(t)$ can exhibit different asymptotic behaviors as $t \to \infty$. 
This term may be constant (e.g., simple shear deformations) or behave as a power-law (e.g.,  dilations in one, two or three dimensions). In analogy with the Boltzmann equation, three distinct asymptotic regimes emerge as $t \to \infty$, depending on the homogeneity $\gamma$ and the type of deformation $L(t)$:

\begin{itemize}
\item The hyperbolic term dominates over the collision term.
\item The collision term dominates over the hyperbolic term.
\item Both terms remain of the same order of magnitude.
\end{itemize}
An outline of the expected behavior for homoenergetic solutions to the Boltzmann equation was presented in \cite{james2019self}. Here, we focus on the Landau equation and provide intuition regarding the long-time behavior of solutions to \eqref{eq:drifted_homo_Landau}. To this end, we rewrite the collision operator in non-divergence form, in components,
\begin{equation}\label{eq:collopL2}
Q[g] = \sum_{i,j}a_{ij}[g]\,\partial_{ij}g + c_\gamma [g]\,g,
\end{equation}
where $a_{ij}[g]$ are the entries of the matrix
$A[g]=(a_{ij}[g])_{i,j=1}^{3}$ (cf.~\eqref{Diffusion-matrixA[g]}), given by 
\begin{equation}
a_{ij}[g](w)=\big(\Phi_{ij}*g\big)(w),
\quad \Phi_{ij}(z)=|z|^{\gamma+2}\Big(\delta_{ij}-\frac{z_iz_j}{|z|^2}\Big),
\end{equation}
with $\Phi$ as in the definition of the collision operator \eqref{eq:collopL} 
and $c_\gamma [g]$ is defined as in \eqref{eq:cgamma}, namely
\begin{equation}
 c_\gamma [g](w)  =
    \begin{cases}
        2(3+\gamma)\big(g*|\cdot|^{\gamma}\big)(w), & \gamma>-3\\
        8\pi g(w), & \gamma = -3.
    \end{cases}
\end{equation}
We now estimate the order of magnitude of the two terms on the right-hand side
of \eqref{eq:collopL2}: 
\begin{align}\label{eq:a_ij}
a_{ij}[g] &\sim \rho(t)\,[w]^{\gamma+2}, \qquad \partial_{ij}g \sim [g]/[w]^2
\quad\Rightarrow\quad a_{ij}[g]\,\partial_{ij}g \sim \rho(t)\,[w]^{\gamma}\,[g]; \\ \label{eq:c}
c_\gamma[g] &\sim \rho(t)\,[w]^{\gamma}
\quad\Rightarrow\quad c_\gamma[g]\,g \sim \rho(t)\,[w]^{\gamma}\,[g].
\end{align}
Since both contributions (cf. \eqref{eq:a_ij}-\eqref{eq:c}) give the same scaling, we obtain that  the collision term scales as  
\begin{equation}
Q(g)(w) \sim \rho(t)\,[w]^{\gamma}\,[g].
\label{eq:scaleQ}
\end{equation}

On the other hand, assuming that $L(t)$ scales as a function of time $[L](t)$, we have that the hyperbolic term rescales as 
\begin{equation}\label{eq:scaleL}
L(t)w\cdot\nabla_{w}g \sim [L](t)[g]
\end{equation} 

Consequently, at a formal level, determining the dominant physical mechanism
reduces to comparing \eqref{eq:scaleQ} and \eqref{eq:scaleL}, namely
$$[L](t) \quad \text{and} \quad \rho(t)[w]^{\gamma}.$$
This formal argument allows one to identify the expected regime
(hyperbolic-dominated, balanced or collision-dominated), which guides the
study of the corresponding long-time behavior. 
We refer to \cite{james2019self} for the rigorous computations in the case of the Boltzmann equation.
Table  \ref{table1} collects the three most important cases, dilation, planar shear and simple shear.

\begin{table}[h]
\centering
\small
\renewcommand{\arraystretch}{1.3}
\begin{tabular}{|l|l|c|l|}
\hline
\textbf{Deformation $L(t)$} & $ \ \mathbf{\rho(t)}$ & $\mathbf{\gamma}$ & \textbf{Regime} \\
\hline
Homogeneous dilation 
 & $\sim t^{-3}$ 
 & $\gamma > -2$ & Hyperbolic-dominated \\
 & & $\gamma = -2$ & Critical \\
 & & $-3 \le \gamma < -2$ & Collision-dominated \\
\hline

Planar shear 
 & $\sim  t^{-1}$ 
 & $\gamma > 0$ & Hyperbolic-dominated \\
 & & $\gamma = 0$ & Balanced (self-similar) \\
 & & $-3 \le \gamma < 0$ & Collision-dominated \\
\hline

Simple shear 
 & const. 
 & $\gamma > 0$ & Collision-dominated \\
 & & $\gamma = 0$ & Balanced (self-similar) \\
 & & $-3 \leq  \gamma < 0$ & Hyperbolic-dominated  \\
\hline
\end{tabular}
\caption{Expected asymptotic regimes for the Landau equation for dilation, planar, and simple shear.}\label{table1}
\end{table}

\FloatBarrier

\medskip 

\bibliography{bibliography}
\bibliographystyle{plain}

\end{document}